\documentclass[12pt]{article}
\usepackage{amssymb,amsmath,amsthm,tikz-cd,amscd,cite,hyperref, pgfplots, tikz, caption, graphicx}
\usepackage[letterpaper, margin=1in]{geometry}

\usepackage{authblk}

\title{The Boundary Harnack Principle and the $3G$ Principle in Fractal-Type Spaces}

\author[1]{Anthony Graves-McCleary} 
\author[2]{Laurent Saloff-Coste}

\affil[1]{Department of Mathematics, Cornell University, Ithaca, NY ag2537@cornell.edu}
\affil[2]{Department of Mathematics, Cornell University, Ithaca, NY lsc@cornell.edu}

\date{}
\newcommand{\set}[1]{\left\{#1\right\}}
\newcommand{\Real}{\mathbf R}

\newcommand{\res}{\textrm{res}}

\newcommand{\N}{\mathbf{N}}

\newcommand{\Z}{\mathbf{Z}}

\newcommand{\ip}[2]{\left< #1, #2  \right>}
\newcommand{\abs}[1]{\left\vert#1\right\vert}
\newcommand{\norm}[1]{\left\Vert#1\right\Vert}

\newcommand{\fct}[3]{#1 \colon #2 \rightarrow #3}

\renewcommand{\phi}{\varphi}

\newcommand{\mc}{\mathcal}

\newcommand{\grad}{\nabla}

\DeclareMathOperator*{\esssup}{ess\,sup}
\DeclareMathOperator*{\essinf}{ess\,inf}

\newtheorem{theorem}{Theorem}[section]
\newtheorem{proposition}[theorem]{Proposition}
\newtheorem{lemma}[theorem]{Lemma}
\newtheorem{definition}[theorem]{Definition}
\newtheorem{assumption}[theorem]{Assumption}
\newtheorem{example}[theorem]{Example}
\newtheorem{corollary}[theorem]{Corollary}

\usetikzlibrary{arrows, matrix}

\numberwithin{equation}{section}

\graphicspath{{C:\Users\antho\OneDrive\Documents\Math Research\Sierpinski_carpet_6.png/}{D:2dCarpet}}

\begin{document}
\maketitle

\begin{abstract}
We prove a generalized version of the $3G$ Principle for Green's functions on bounded inner uniform domains in a wide class of Dirichlet spaces. In particular, our results apply to higher-dimensional fractals such as Sierpinski carpets in $\Real^n$, $n\geq 3$, as well as generalized fractal-type spaces that do not have a well-defined Hausdorff dimension or walk dimension. This yields new instances of the $3G$ Principle for these spaces. We also discuss applications to Schr\"odinger operators.
\end{abstract}
\section{Introduction}

%maybe don't need to spend a page on Dirichlet forms. Reference something (Mathav). 

%Applications section: mention Sierpinski gasket and recover Hansen's result using our main theorem. Then point out Mathav's example (Kumagai, Barlow, Coulon paper on cable trees)

%Given an abstract domain $\Omega$ in a metric space and a function $G$ on $\Omega\times \Omega$ with properties mimicking a Green's function, we are interested in the equivalence between the Uniform Boundary Harnack Principle for $G$ and the $3G$ inequality for $G$, also known as the generalized triangle property. Hansen discussed this in the setting when $G(x,y)$ resembles $d(x,y)^{-\gamma}$, where $d$ is the metric. This includes the classical Euclidean Green's function in inner domains as well as the Green's function for fractional Laplacian.\\

%We generalize the expression $d(x,y)^{-\gamma}$ in two ways. First, non-classical exponents $\gamma$ can be further generalized using a scale function (see Definition \ref{scaledef}). Next, while Hansen's paper \cite{hansenmain} does not introduce a notion of volume in the abstract setting, the expression $d(x,y)^{-\gamma}$ for a Green's function in a Dirichlet space implicitly assumes a certain amount of uniform behavior of volume of balls. We can avoid this by introducing a measure to the abstract setting and express $G$ in terms of it and the scale function.\\

We are interested in the behavior of a Green's function $G(x,y)$ in a domain at points near the boundary of the domain, as well as ratios of the Green's function across the entire domain. As a warm up and to introduce many of the concepts involved, we start with a section in Euclidean space. Then we make the leap to the abstract setting, working in domains in metric measure spaces with no \textit{a priori} notion of harmonicity. It is here that we prove our main theorem \ref{maithm}. Then in the following sections we introduce concepts that will yield many examples: inner uniform domains and Dirichlet spaces. The paper by Murugan \cite{murugan1} gives a rich collection of examples for which the Generalized $3G$ Principle, a concept linked to the main theorem and defined in \ref{3g} below, was not previously known to hold.

\subsection{Notation}

Let $a\vee b=\min(a,b)$ and $a\wedge b=\max(a,b)$. We let $a\preceq b$ denote the existence of a constant $C>0$ such that $a\leq C\cdot b$. We let $a\asymp b$ denote that both $a\preceq b$ and $a\succeq b$ hold. Exact constants will sometimes change line to line. On occasion we will write $a:b$ to denote $a/b$. Given a metric space $X$, we denote by $C_c(X)$ the set of all continuous, compactly supported real-valued functions on $X$. Given $M>0$ and a ball $B=B(x,r)$ of center $x$ and radius $r>0$, $MB$ will denote the ball $B(x, Mr)$ of center $x$ and radius $Mr$.

\section{Classical Preliminaries and Motivation}

For this section, work in the open unit ball $B$ in $\Real^n$ with $n\geq 3$. Let $\alpha(n)$ denote the volume of the unit ball in $\Real^n$ and let $c_n=\frac{1}{n(n-2)\alpha(n)}$. For $x,y\in B$, $x\neq y$ and $y\neq 0$ let $$G_B(x,y)=c_n\left(\norm{x-y}^{2-n} - \left(\norm{y}\norm{x-\frac{y}{\norm{y}^2}}\right)^{2-n}\right).$$ Then if $x\neq 0$ let $$G_B(x,0)=c_n\left(\norm{x}^{2-n}-1\right).$$ This is the classical Green's function for $B$. Note that for any $\xi\in \partial B$ and any $x\in B$ we have $\lim_{y\rightarrow \xi} G_B(x,y)=0$. In fact, for $x\in B$ fixed and $y\in B$ sufficiently close to $\partial B$, $G_B(x,y)$ resembles the distance to the boundary $d(y,\partial B)$, up to constants.\\

In other domains such as cones, the Green's function $G(x,y)$ may have different boundary decay, but in classical settings it can controlled by the Uniform Boundary Harnack Principle. Let $D$ be a domain in $\Real^n$. With an eye towards generalizing, write $d(x,y)=\norm{x-y}$. 

\begin{definition} Let $D$ be a domain let $G$ be a function on $D\times D$. We say that $D$ satisfies the \textbf{Uniform Boundary Harnack Principle} for $G$ if there exist constants $r_0, M, C>0$ such that for all $\xi\in \partial D$, $0<r\leq r_0$, and for all $x_1, x_2, y_1, y_2\in D$ with $Md(x_j, \xi)<r\leq d(y_j, \xi)$, \begin{equation}\label{boundaryineq}\frac{G(x_1, y_1)}{G(x_2, y_1)}\leq C\frac{G(x_1, y_2)}{G(x_2, y_2)}.\end{equation}
\end{definition}

Note that in the inequality \ref{boundaryineq}, $y_1$ and $y_2$ can be switched. The definition of the Uniform Boundary Harnack Principle for more abstract domains is exactly the same. This is all stated in terms of $G$, but when $G$ is truly the Green's function of the domain $D$, this implies a similar bound for any two positive harmonic functions $u,v>0$ defined around $\xi$ that vanish on $\partial D$. See Lierl \cite{lierl1} for further details.\\

The Euclidean open unit ball $B$ satisfies the Uniform Boundary Harnack Principle for its Green's function $G_B$. The Uniform Boundary Harnack Principle has a rich history and has been generalized with regard to both the domain and the kind of Green's function it concerns. It originated in papers by Ancona \cite{ancona1}, Dahlberg \cite{dahlberg1} and Wu \cite{wu1} in the late 1970s. See Armitage and Gardiner \cite{classicalpotential} section 8.7 for a textbook presentation in a Lipschitz domain in Euclidean space. See Barlow and Murugan \cite{barlowmurugan1} for a more comprehensive review as well as a state-of-the-art proof on inner uniform domains in Dirichlet spaces satisfying only an elliptic Harnack inequality.\\

Our other central concept, the $3G$ Principle, arises from the desire to control ratios of the Green's function on the entire domain. This is particularly important in the study of conditional Brownian motion and Schr\"odinger operators, which we will discuss in section 6. In the ball, there is a constant such that for any $x,y,z\in B$, \begin{equation}\label{3gclassical} \frac{G_B(x,z)G_B(y,z)}{G_B(x,y)} \leq C\left(\norm{x-z}^{2-n}+\norm{y-z}^{2-n}\right).\end{equation}

Notably $\norm{x-z}^{2-n}$ is, up to a constant, the Green's function of all of $\Real^n$. Therefore this ``$3G$" ratio of Green's functions in the ball can be controlled by Green's functions from a larger domain. In generalizing this, we can avoid reference to any such larger domains by use of an additional benchmark Green's function. We are now ready to state our version of the $3G$ Principle.

\begin{definition}\label{3g} Let $D$ be a domain and let $G$ be a function on $D\times D$. Fix $o\in D$ and let $g(x)= G(x,o)\wedge 1$ be the benchmark Green's function. We say that $D$ satisfies the \textbf{Generalized $3G$ Principle} with respect to $G$ if there exists a constant $C>0$ such that for all $x,y,z\in D$ we have \begin{equation}
\frac{G(x,z) G(y,z)}{G(x,y)} \leq C\left(\frac{g(z)}{g(x)}G(x,z)+\frac{g(z)}{g(y)}G(z,y)\right).
\end{equation}
\end{definition}

In practice this will be independent of choice of $o\in D$, up to constants. The presence of $g$ is necessary in many cases due to boundary vanishing of the Green's function. In applications, the expression $\frac{g(z)}{g(x)}G(x,z)$ will be bounded above by the Green's function of the entire space, therefore relating this back to \ref{3gclassical}. See Lemma \ref{global3g} below for more details. The \ref{3gclassical} form of the $3G$ Principle appears for bounded Lipschitz domains in Euclidean space in Cranson, Fabes, and Zhao \cite{cranston1} with a view towards Schr\"odinger operators, which we will discuss in a later section. The Generalized $3G$ Principle appears in Riahi \cite{riahi1}. The proofs of the $3G$ Principle almost always make use of the Uniform Boundary Harnack Principle, which raises the question of whether they are equivalent. Hansen \cite{hansenmain} demonstrated the equivalence between the Uniform Boundary Harnack Principle and a form of the $3G$ Principle (see Definition \ref{strong3g} below) in an abstract setting that, in particular, applies to bounded uniform domains in $\Real^n$, $n\geq 3$.\\

 We now state an alternate form of the Generalized $3G$ Principle which is routine to verify as equivalent to the original.

\begin{proposition}\label{alt3g} Let $D$ be a domain and let $G$ be a function on $D\times D$. Fix $o\in D$ and let $g(x)=G(x,o)\wedge 1$ be the benchmark Green's function. Also let $\widetilde{G}(x,y)=G(x,y)/(g(x)g(y))$. Then $D$ satisfies the Generalized $3G$ Principle for $G$ with constant $C>0$ if and only if for all $x,y,z\in D$, \begin{equation}\frac{1}{\widetilde{G}(x,y)}\leq C\left(\frac{1}{\widetilde{G}(x,z)}+\frac{1}{\widetilde{G}(z,y)}\right).\end{equation}

\end{proposition}

\noindent\textbf{Proof:} Trivial.\qed\\

This alteration is the same as stating that $\widetilde{G}(x,z)\wedge \widetilde{G}(z,y)\leq C\widetilde{G}(x,y)$. This suggests the following strengthening:

\begin{definition}\label{strong3g} Let $D$ be a domain and let $G$ be a function on $D\times D$. Fix $o\in D$ and let $g(x)=G(x,o)\wedge 1$ be the benchmark Green's function. We say that $D$ satisfies the \textbf{Strong Generalized $3G$ Principle} with respect to $G$ if there exists a constant $C>0$ such that for all $x,y,z\in D$ such that $d(z,x)\leq d(z,y)$ we have \begin{equation}
\widetilde{G}(z,y)\leq C\widetilde{G}(x,y).
\end{equation}
\end{definition}

It is this version of the $3G$ Principle that we will focus on in our main theorem \ref{maithm}.

%\subsection{Finding the Correct Domains to Study}

To close this section, let us remark about the kinds of domains we will study. In Euclidean space, the Uniform Boundary Harnack Principle and the Generalized $3G$ Principle do not hold on all domains even in Euclidean space with classical Green's function. We will impose some structure on our domains to get results, although it is a worthwhile question whether our assumptions can be weakened. Up this point, the Uniform Boundary Harnack Principle has mainly been studied on inner uniform domains, to be defined in section 5. Some Euclidean domains, such as the slit disk $D=\set{(x,y)\in \Real^2\colon x^2+y^2<1}\setminus\set{(x,0)\colon 0\leq x<1}$ are inner uniform with a boundary different from their usual Euclidean boundary.\\

Inner uniformity requires the notion of length of curves. We avoid this for now by substituting the interior corkscrew condition plus the (classical) Harnack inequality, defined later in \ref{bigassumption}.

\section{Abstract Setting}

In this section we establish the abstract setting for our main theorem \ref{maithm}. We are working with the least amount of information needed to express a Green's function as well as control its behavior. Notably we will not make direct reference to harmonicity, so we eschew the use of Dirichlet spaces in this section. They will appear in the examples. Our ambient space will be a metric measure space:

\begin{definition} 

A \textbf{metric measure space} is a triple $(X, d, \mu)$ where $(X, d)$ is a metric space and $\mu$ is a Borel measure on $X$.

\end{definition}

Given a metric measure space $(X, d, \mu)$, let $B(x,r)=\set{y\colon d(x,y)<r}$ and let $V(x,r)=\mu(B(x,r))$. Given a subset $D\subseteq X$ we let $\partial D$ denote its boundary and $\delta(\cdot)=d(\cdot, \partial D)$. Given $y\in D$ we let $B(y)=B(y, \delta(y)/2)$ and $D(y)=D\setminus B(y)$.\\

To work with our abstract Green's function, we need some amount of control over the volume of balls.

\begin{definition}
 A metric measure space $(X, d, \mu)$ satisfies \textbf{volume doubling (VD)} if there exists a constant $C>0$ such that for all $x\in X$ and $r>0$, $$V(x, 2r)\leq CV(x,r).$$

\end{definition}

\begin{definition}\label{scaledef}

A \textbf{scale function} for a metric measure space $(X, d, \mu)$ is a continuous increasing bijection $\fct{\Psi}{(0, \infty)}{(0, \infty)}$ for which there exist constants $1<\beta\leq \beta'$ and $C>1$ such that for all $0<r\leq R$, we have \begin{equation}C^{-1}\left(\frac{R}{r}\right)^\beta\leq \frac{\Psi(R)}{\Psi(r)}\leq C\left(\frac{R}{r}\right)^{\beta'}.\end{equation}

\end{definition}

\begin{figure}
%\begin{center}
\centering
\begin{tikzpicture}

\draw (1.5,2) -- (4,2) -- (3,-1) -- (-1,-1) -- (0,-3) -- (-4, -3) -- (-7, 4) -- (2,4) -- (1.5,2);
\draw (2.5, 0) node{$D$};
\draw (-2, 2) node{$\bullet$};
\draw (-1.75,2.1) node{$y$};
\draw (-2.3, 1.44) node{$B(y)$};
\draw (-4, 1.44) node{$D(y)$};
\draw (-2, 2) circle [radius=1cm];
\draw (-2, 2) -- (-2,4);

\end{tikzpicture}
\captionof{figure}{A set $D$ and $y\in D$ with $B(y)$ and $D(y)$ labeled. The vertical line segment has length $\delta(y)$.}
%\end{center}
\end{figure}
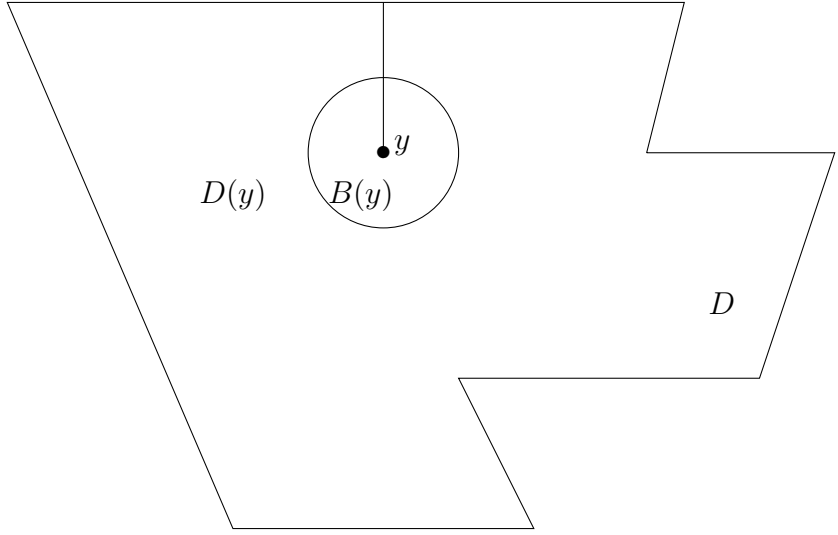

The scale and volume functions inspire the following general definition.

\begin{definition}
A function $\fct{f}{[0,\infty)}{[0,\infty)}$ is \textbf{doubling} if there exists a constant such that $f(2r)\leq Cf(r)$ for all $r\geq 0$.
\end{definition}

For example, if $(X, d, \mu)$ satisfies (VD), then $r\mapsto V(x,r)$ is doubling for each $x\in X$, with a constant that is uniform in $x\in X$. If $\Psi$ is a scale function, then $\Psi$ is doubling. This doubling property will be of vital importance.\\

In this section let $(X, d,\mu)$ be a metric measure space satisfying (VD). Let $\Psi$ be a scale function for $(X, d,\mu)$. Let $D\subseteq X$ be open and precompact. We will consider a function $\fct{G}{D\times D\setminus \set{(x,x)\colon x\in D}}{(0, \infty)}$ along with the following assumptions on $D$ and on $G$:

\begin{assumption}\label{bigassumption} There exist constants $M\geq 3$, $r_0>0$, $C_0>0$, $\gamma\geq 0$ such that the following four properties hold:\\

\textbf{(i)} Interior Corkscrew Property: For every $\xi\in \partial D$ and all $0<r\leq r_0$, there exists $x\in D$ such that $d(x,\xi)<r$ and $\delta(x)>r/M$.\\

\textbf{(ii)} Quasi-symmetry: For all $x,y\in D$, \begin{equation}G(x,y)\leq C_0 G(y,x).\end{equation}

\textbf{(iii)} Quasi-polynomial decay: For all $x,y\in D$, $x\neq y$, \begin{equation}G(x,y)\leq C_0\frac{\Psi(d(x,y))}{V(x,d(x,y))}.\end{equation} Also, for $x\in \overline{B(y)}$, \begin{equation}G(x,y)\geq \frac{1}{C_0}\frac{\Psi(d(x,y))}{V(x,d(x,y))}.\end{equation}

\textbf{(iv)} Harnack's inequality: For all $y\in D$ and $x_1, x_2\in D(y),$ \begin{equation}G(x_1, y)\leq C(k)G(x_2, y)\end{equation} when $d(x_1, x_2)\leq k \delta(x_1)\wedge \delta(x_2),$ $k\in \N$.\\

\textbf{(v)} For all $x,y,z\in D$, all distinct, $d(x,y)\leq d(y,z)$, we have \begin{equation}
\frac{\Psi(d(y,z))}{\Psi(d(x,y))}\leq C\frac{V(z, d(y, z))}{V(x, d(x,y))}.\end{equation} 
\end{assumption}

Note that if \ref{bigassumption} holds for $(X, d, \mu)$, $\Psi$, and $G$, then \ref{bigassumption} also holds for $(X, d', \mu)$, $\Psi$ and $G$ where $d'=d/C'$, $C'>0$ a constant, albeit with different constants $M, R, r_0, C_0$. This is a consequence of the volume doubling property and properties of the scale function. Therefore when working with \ref{bigassumption} we may divide $d$ by the diameter of $D$ to assume that $D$ has metric diameter $1$. We are ready to state our main theorem.

\begin{theorem}\label{maithm} Let $(X, d,\mu)$, $\Psi$, $D$ and $G$ satisfy Assumption \ref{bigassumption}. Then $D$ satisfies the Uniform Boundary Harnack Principle with respect to $G$ if and only if $D$ satisfies the Strong Generalized $3G$ Principle with respect to $G$. 
\end{theorem}

Note that implicitly in the statement of the Strong Generalized $3G$ Principle we have picked a distinguished point $o\in D$. As one would expect, the proof will require extreme care in studying points near the boundary of $D$. We introduce additional concepts in this direction. Given $x,y\in D$, let \begin{equation}r(x,y)= \delta(x)\vee \delta(y)\vee d(x,y).\end{equation} The function $r(x,y)$ is a way of simultaneously tracking the three most important metric quantities related to $x$ and $y$: their distances to the boundary and to each other. Next, we need a scale of smallness. We choose \begin{equation} \epsilon:= \frac{r_0}{12M}.\end{equation} The exact value of $\epsilon$ is not important, but we need it to be smaller than $r_0$. Next, given $x,y\in D$ we define \begin{equation}\mathcal{B}(x,y)=\set{A\in D\colon \delta(A)>\frac{r(x,y)}{M}, d(x,A)\vee d(y,A)<5r(x,y)}\end{equation} if $r(x,y)<\epsilon$, and $\mathcal{B}(x,y)=\set{o}$ otherwise.\\

To motivate this definition, $r(x,y)<\epsilon$ will turn out to be a quite delicate situation, thus requiring us to make this distinction. The points $A\in \mathcal{B}(x,y)$ are ``close to $x$ and $y$, far from $\partial D$" on the scale of $r(x,y)$.\\

In view of the Interior Corkscrew Property \ref{bigassumption}(i), let us define for $\xi\in \partial D$ and $0<\rho\leq r_0$, \begin{equation}\label{corkset} \mathcal{A}_\rho(\xi)=\set{A\in D\colon d(A, \xi)<\rho, \delta(A)>\rho/M}.\end{equation}

Note that $\mathcal{A}_\rho(\xi)$ is always nonempty by \ref{bigassumption}(i). We have the following relation with the sets $\mathcal{B}(x,y)$.

\begin{lemma}\label{cones}
Let $x,y\in D$ with $r(x,y)<\epsilon$ and let $\xi,\eta\in \partial D$ such that $d(x,\xi)=\delta(x),$ $d(y, \eta)=\delta(y)$. Then \begin{equation}
\mathcal{A}_{r(x,y)}(\xi)\cup \mathcal{A}_{r(x,y)}(\eta)\subseteq \mathcal{B}(x,y).
\end{equation}
\end{lemma}

\noindent\textbf{Proof:} If $A\in \mathcal{A}_{r(x,y)}(\xi)$ then $\delta(A)>r(x,y)/M$, so $$d(x,A)\leq d(x,\xi)+d(\xi, A)<2r(x,y)$$ and $$d(y,A)\leq d(y,x)+d(x,A)<3r(x,y).$$ Similarly for $\eta$.\qed\\

As before, for a distinguished point $o\in D$ we define $g(x)=1\wedge G(x,o)$. We also record the following basic fact.

\begin{lemma}
Let $x,y\in D$ and let $A_1, A_2\in \mc{B}(x,y)$. Then \begin{equation} g(A_1)\asymp g(A_2)\end{equation} with constants independent of $x,y$.

\end{lemma}

\noindent\textbf{Proof:} Let $x,y\in D$ and $A_1, A_2\in \mc{B}(x,y)$. If $r(x,y)\geq \epsilon$, then $A_1=A_2=o$ so $g(A_1)=g(A_2)$. If $r(x,y)<\epsilon$, then $d(A_1, A_2)\leq d(x,A_1)+d(x,A_2)\leq 10r(x,y)\leq 10M \delta(A_1)\wedge \delta(A_1)$ and the result follows from Harnack's inequality \ref{bigassumption}(iv).\qed\\

In proving the main theorem \ref{maithm}, it will be useful to introduce an intermediate property, the $\mathcal{B}$-approximation property, that will be proved equivalent to both the Uniform Boundary Harnack Principle and the Generalized $3G$ Principle.

\begin{definition} Under Assumption \ref{bigassumption}, we say that the \textbf{$\mathcal{B}$-approximation property} holds if there exist constants $c,C>0$ such that for all $x,y\in D$ and every $A\in \mathcal{B}(x,y)$, \begin{equation} c\frac{\Psi(d(x,y))}{V(x,d(x,y))g(A)^2}\leq \frac{G(x,y)}{g(x)g(y)}\leq C\frac{\Psi(d(x,y))}{V(x,d(x,y))g(A)^2}.\end{equation}
\end{definition}

As we just stated above, we have the following theorems. They are intermediate steps in proving the main theorem \label{mainthm}.

\begin{theorem}\label{mainpart2}
Under Assumption \ref{bigassumption}, the Uniform Boundary Harnack Principle holds if and only if $\mathcal{B}$-approximation holds.
\end{theorem}

\begin{theorem}\label{mainpart3}
Under Assumption \ref{bigassumption}, the Generalized $3G$ Principle holds if and only if $\mathcal{B}$-approximation holds.
\end{theorem}

Let us finish this section with some important preparatory lemmas. All statements are under Assumption \ref{bigassumption} unless otherwise stated. 

\begin{lemma}\label{amazinglemma} $g\asymp 1$ on $\set{\delta \geq \epsilon/(8M^3)}$.
\end{lemma}

\noindent\textbf{Proof:} By Assumption \ref{bigassumption} (iii), for $x\in \overline{B(o)}$ we have \begin{equation} G(x,o)\geq \frac{1}{C_0}\frac{\Psi(d(x,o))}{V(x,d(x,o))} \asymp \frac{\Psi(d(x,o))}{V(o,d(x,o))},\end{equation} where the last estimate follows from volume doubling. Thus $g\asymp 1$ on $B(o)$.\\

Next, fix $z\in \partial B(o)$. (This exists by connectedness of $D$.) Then $d(\cdot, x)\leq 1\leq (8M^3/\epsilon)\delta$ on $K=\set{\delta\geq \epsilon/(8M^3)}$. By Harnack's inequality \ref{bigassumption}(iv), $G(\cdot, o)\asymp G(z,o)=g(z)\asymp 1$ on $K\setminus B(o)$ whence $g\asymp 1$ on $K\setminus B(o)$ as well.\qed\\

\begin{lemma}\label{hansenl2.2} Let $x,y\in D$ such that $d(x,y)\leq (8M^3/\epsilon)\delta(x)\wedge \delta(y)$. Then \begin{equation}G(x,y)\asymp \frac{\Psi(d(x,y))}{V(x, d(x,y))}.\end{equation}\end{lemma}

\noindent\textbf{Proof:} The upper bound follows from \ref{bigassumption}(iii), so we focus on the lower bound.\\

(a) If $d(x,y)\leq \delta(y)/2$ then the lower bound from \ref{bigassumption}(iii) gives \begin{equation}G(x,y)\succeq \frac{\Psi(d(x,y))}{V(x,d(x,y))}.\end{equation}

(b) If $d(x,y)>\delta(y)/2$, so consequently $d(x,y)\asymp \delta(y)$. Choose a point $z\in \partial B(y, \delta(y)/2)$. Then \begin{equation}d(x,z)\leq d(x,y)+\frac{\delta(y)}{2}<2d(x,y)\leq \frac{16 M^3}{\epsilon} \delta(x),\end{equation} where $\delta(x)\leq \delta(y)\leq 2\delta(z).$ By Harnack's inequality \ref{bigassumption}(iv) and (a), we conclude that \begin{equation}G(x,y)\asymp G(z,y)\asymp \frac{\Psi(d(z,y))}{V(z, d(z,y))} \asymp \frac{\Psi(d(x,y))}{V(y,d(x,y))},\end{equation} the last following from the fact that $d(z,y)=\delta(y)/2\asymp d(x,y)$ and the doubling property for $V$ as well as the doubling property of $\Psi$.\qed\\

\section{Proof of the Main Theorem}

We begin with the proof that the $\mathcal{B}$-approximation property implies the Uniform Boundary Harnack Principle.\\

\noindent\textbf{Proof of Theorem \ref{mainpart2}, $\Leftarrow$ direction:} Assume \ref{bigassumption} and assume that $D$ satisfies $\mc{B}$-approximation.\\

 Fix $\xi \in \partial D$, $0<r\leq r_0$, and $x_1, x_2, y\in D$ with $Md(x_j, \xi)<r\leq d(y, \xi)$. By the $\mc{B}$-approximation property, for every choice of $A_j\in \mc{B}(x_j, y)$, \begin{equation}\frac{G(x_1, y)}{G(x_2, y)}\asymp \frac{g(x_1)g(y) \Psi(d(x_1,y))}{g(A_1)^2 V(x_1,d(x_1,y))}\colon \frac{g(x_2)g(y)\Psi(d(x_2, y))}{g(A_2)^2 V(x_2, d(x_2, y))}.\end{equation}

Now $r(1-\frac{1}{M})\leq d(x_j, y)\leq 1$ so that $d(x_1, y)\asymp d(x_2, y)$. Further, $d(x_1, x_2)\leq 2r/M\preceq d(x_j, y).$ Therefore by volume doubling and doubling for $\Psi$ we have $V(x_1, d(x_1, y))\asymp V(x_2, d(x_2, y))$. Thus we have that \begin{equation}\frac{g(x_1)g(y) \Psi(d(x_1,y))}{g(A_1)^2 V(x_1,d(x_1,y))}\colon \frac{g(x_2)g(y)\Psi(d(x_2, y))}{g(A_2)^2 V(x_2, d(x_2, y))} \asymp \frac{g(x_1)}{g(x_2)}\colon \frac{g(A_1)^2}{g(A_2)^2}.\end{equation}\\

Thus to prove the Uniform Boundary Harnack Principle, it suffices to produce a uniform bound on $g(A_1)/g(A_2)$.\\

If $r(x_1, y)\geq \epsilon/2$ and $r(x_2, y)\geq \epsilon/2$, then $\delta(A_j)>\epsilon/(2M)$ for $j=1, 2$ and therefore $g(A_1)\asymp 1 \asymp g(A_2)$ by Lemma \ref{amazinglemma}.\\

So let us suppose that \begin{equation}\rho:= 2 r(x_1, y)\wedge r(x_2, y)<\epsilon,\end{equation} and choose $A\in \mc{A}_\rho(\xi)$. We intend to show that $A\in \mc{B}(x_1,y)\cap\mc{B}(x_2, y)$ (and this will finish the proof). To that end we may assume without loss of generality that \begin{equation}r(x_1, y)\leq r(x_2, y).\end{equation}

Let us recall that, by our assumption at the beginning of this proof, $Md(x_j, \xi)<r<d(y, \xi)$. This implies that, for $j=1,2$, \begin{equation}r(x_j, y)\geq d(x_j, y)\geq d(y, \xi)-d(x_j, \xi)>r-\frac{r}{M}\geq \frac{2r}{M}\end{equation} whence \begin{equation}d(x_j, A)\leq d(x_j, \xi)+d(\xi, A)<\frac{r}{M}+\rho<3r(x_j, y),\end{equation} $$d(y,A)\leq d(y, x_j)+d(x_j, A)<4r(x_j, y).$$

Indeed, this is because $d(x_j, \xi)=\delta(x_j)\leq r(x_j, y)$, $\rho=2 r(x_1, y)\wedge r(x_2, y)\leq 2r(x_j, y)$. Also, $d(y, x_j)\leq r(x_j, y)$ and the rest follows.\\

Furthermore we have \begin{equation}\delta(A)>\frac{\rho}{M}= \frac{2r(x_1, y)}{M}.\end{equation} So certainly $A\in \mc{B}(x_1, y)$. To prove that $A\in \mc{B}(x_2, y)$ it remains to show that $2r(x_1, y)\geq r(x_2, y)$. Obviously, \begin{equation}\delta(x_2)\vee d(x_1, x_2)\leq d(x_1, \xi)+d(x_2, \xi)<\frac{2r}{M}<r(x_1, y),\end{equation} where the last inequality follows from earlier estimates on $r(x_j, y)$. Consequently, $d(x_2, y)\leq d(x_2, x_1)+d(x_1, y)<2r(x_1, y)$. Since $\delta(y)\leq r(x_1, y)$, we finally conclude that \begin{equation}r(x_2, y)=d(x_2, y)\vee \delta(x_2)\vee\delta(y)<2r(x_1, y).\end{equation}

Therefore $A\in \mc{B}(x_2, y)$. We can therefore use $\mc{B}$-approximation to conclude that the uniform boundary Harnack inequality holds. This concludes the proof of the $\Leftarrow$ direction of Theorem \ref{mainpart2}.\qed\\

\noindent\textbf{Proof of Theorem \ref{mainpart2}, $\Rightarrow$ direction:} Let us assume that $D$ satisfies the Uniform Boundary Harnack Principle with respect to $G$. Our goal for the $\Rightarrow$ direction is to prove that $\mc{B}$-approximation holds. This is the harder direction and will be accomplished in parts, according to the positioning of points $x, y\in D$ with respect to the boundary and to each other.\\

Let us first observe the following: If $x\in \set{\delta<6\epsilon}$, then $d(x,o)\geq \delta(o)-6\epsilon \geq \delta(o)/2$. Therefore \begin{equation}\set{\delta<6\epsilon}\subseteq D(o),\end{equation} \begin{equation}g=G(\cdot, o) \textrm{ on } \set{\delta<6\epsilon}.\end{equation} This will be used without being mentioned explicitly.\\

Fix $x,y\in D$ and $A\in \mc{B}(x,y)$.\\

\noindent\textbf{Claim:} \begin{equation}\label{claim1} g(A)\asymp 1 \textrm{ when } r(x,y)\geq \frac{\epsilon}{8M^2}.\end{equation}

\noindent\textbf{Proof of Claim:} Indeed, if $r(x,y)\geq \epsilon$, then $A=o$, $g(A)=1$. Suppose that $r(x,y)<\epsilon$. Then $\delta(A)>r(x,y)/M\geq \epsilon/(8M^3)$ so this follows from Lemma \ref{amazinglemma}.\qed\\

Without loss of generality, assume \begin{equation}\delta(x)\leq \delta(y).\end{equation}

Let us choose points $\xi, \eta\in \partial D$ such that $d(x, \xi)=\delta(x)$ and $d(y, \eta)=\delta(y)$, respectively. Several cases will be considered:\\

\noindent\textbf{(a)} (Both $x$ and $y$ far from $\partial D$) Suppose for this case that $\delta(x)\geq \epsilon/(8M^2)$. Then $g(x)\asymp g(y)\asymp 1 \asymp g(A)$ be Lemma \ref{amazinglemma}. Moreover, $G(x,y)\asymp \frac{\Psi(d(x,y))}{V(x,d(x,y))}$ by Lemma \ref{hansenl2.2}. Thus $\mc{B}$-approximation holds in this case.\\

\noindent\textbf{(b)} ($x$ close to $\partial D$ and $y$ far) Suppose for this case that $\delta(x)<\epsilon/(8M^2)$ and $\delta(y)\geq \epsilon/(2M)$. Let $r=\epsilon/(8M)$ and choose $x_1\in \mc{A}_{r/M}(\xi)$. By Lemma \ref{amazinglemma} and by the claim \ref{claim1}, \begin{equation} g(x_1)\asymp g(y)\asymp 1\asymp g(A).\end{equation}

Obviously $d(x_1, \xi)<r/M$, $d(x,\xi)=\delta(x)<r/M$, and $r<\delta(y)\leq d(y,\xi)$. Thus, by the Uniform Boundary Harnack principle, \begin{equation}\frac{G(x,y)}{g(x)}\asymp \frac{G(x_1, y)}{g(x_1)},\end{equation} where $G(x_1, y)\asymp \frac{\Psi(d(x_1, y))}{V(x_1, d(x_1, y))}$ by Lemma \ref{hansenl2.2}. Moreover, \begin{equation}\frac{d(x,y)}{2}-d(x, x_1)\geq \frac{\delta(y)-\delta(x)}{2}-\left(\delta(x)+\frac{r}{M}\right)\geq \frac{\epsilon}{M}\left(\frac{1}{4}-\frac{3}{2}\cdot \frac{1}{8M}-\frac{1}{8M}\right)>0.\end{equation}

Therefore $d(x_1, y)\asymp d(x,y)$: indeed, $d(x_1, y)\leq d(x_1, y)+d(x,y)\leq 2d(x,y)+d(x,x_1)\leq \frac{5}{2}d(x,y)$, and conversely $d(x,y)\leq d(x_1, y)+d(x_1, x)\leq d(x_1, y)+\frac{d(x,y)}{2}$ so that $d(x,y)\leq d(x_1, y)$. The $\mc{B}$-approximation property follows for case (b).\\

\noindent\textbf{(c)} ($x$ and $y$ both potentially close to $\partial D$, but close to each other) Assume for this case that $\delta(y)<\epsilon$ and $d(x,y)<\delta(y)/2$, i.e. $x\in B(y)$. Then \begin{equation}G(x,y)\asymp \frac{\Psi(d(x,y))}{V(x,d(x,y))}.\end{equation} Moreover, $\delta(x)\geq \delta(y)/2$ and $r(x,y)=\delta(y)<\epsilon$ whence $\delta(A)>r(x,y)/M$ and \begin{equation}d(y,A)\leq 5r(x,y)\leq 5M \delta(y)\wedge \delta(A).\end{equation} In addition, \begin{equation}\delta(A)\leq \delta(y)+d(y,A)<\delta(y)+5r(x,y)<6\epsilon.\end{equation} Therefore by the Harnack inequality we find that \begin{equation}g(x)\asymp g(y)\asymp g(A).\end{equation} This proves the $\mc{B}$-approximation in case (c).\\

\noindent\textbf{(d)} Having established (a)-(c), we may assume in the following that $$\delta(x)<\frac{\epsilon}{8M^2},$$ $$\delta(x)\leq \delta(y)<\frac{\epsilon}{2M},$$ $$d(x,y)>\frac{\delta(y)}{2}.$$

\noindent\textbf{(e)} ($x$ and $y$ far apart, but not overly far) Suppose for this case that $d(x,y)\leq 6M^2 \delta(x)$. Then \begin{equation}d(x,y)\asymp \delta(y)\asymp \delta(x)\asymp r(x,y)\end{equation} and $\delta(y)\leq d(x,y)+\delta(x)\leq (6M^2+1)\delta(x)<\epsilon$ whence $r(x,y)<\epsilon$. Therefore $\delta(A)\asymp r(x,y)$ and $d(x,A)<5r(x,y)$. So we have \begin{equation}g(y)\asymp g(x)\asymp g(A).\end{equation} Fix $z\in \partial B(y, \delta(y)/2)$. Then $\delta(z)\geq \delta(y)/2\geq \delta(x)/2$, $d(x,z)\leq d(x,y)+\delta(y)/2\leq 2d(x,y)$, and therefore $d(x,z)\leq 24M^2\delta(x)\wedge \delta(z)$. By \ref{bigassumption}(iii) and (iv) as well as (VD) and ($\Psi$D), \begin{equation}G(x,y)\asymp G(z,y)\asymp \frac{\Psi(d(z,y))}{V(z, d(z,y))}\asymp \frac{\Psi(\delta(y))}{V(z, \delta(y)/2)}\asymp \frac{\Psi(d(x,y))}{V(x, d(x,y))}.\end{equation} Thus $\mc{B}$-approximation holds in the case (e).\\

\noindent\textbf{(f)} For this case let $6M^2\delta(x)<d(x,y)\leq 2M\delta(y)$. Then $r(x,y)\leq 2M\delta(y)<\epsilon$, \begin{equation}d(x,y)\asymp \delta(y)\asymp r(x,y)\asymp \delta(A).\end{equation} Taking $r=\delta(y)/3$ we have \begin{equation}d(x, \xi)=\delta(x)<\frac{\delta(y)}{3M}=\frac{r}{M}\end{equation} and \begin{equation}d(y, \xi)-r\geq d(x,y)-r-\delta(x)>\left(\frac{1}{2}-\frac{1}{3}-\frac{1}{3M}\right)\delta(y) \geq 0.\end{equation}\\

Fix $x_1\in \mc{A}_{r/M}(\xi)$. Then $\delta(x_1)>r/M^2$ and $d(x_1, \xi)<r/M$. By the Uniform Boundary Harnack Principle, \begin{equation}\frac{G(x,y)}{g(x)}\asymp \frac{G(x_1, y)}{g(x_1)}.\end{equation}

Moreover, using $\delta(y)< 2d(x,y)$, \begin{equation}\abs{d(x,y)-d(x_1, y)}\leq d(x, x_1)<\delta(x)+\frac{r}{M}<\left(\frac{1}{6M^2}+\frac{2}{3M}\right)d(x,y)<\frac{d(x,y)}{3}\end{equation} whence $\frac{2}{3}d(x,y)<d(x_1, y)<\frac{4}{3}d(x,y).$ Furthermore, $d(x,y)\leq 2M \delta(y) \leq 6M^3 \delta(x_1).$ So, by Lemma \ref{hansenl2.2} as well as (VD) and ($\Psi$D), \begin{equation}G(x_1, y)\asymp \frac{\Psi(d(x_1, y))}{V(x_1, d(x_1, y))}\asymp \frac{\Psi(d(x, y))}{V(x, d(x,y))}.\end{equation}\\

Moreover, Harnack's inequality \ref{bigassumption}(iv) yields that \begin{equation}g(y)\asymp g(A)\asymp g(x_1).\end{equation} Thus \begin{equation}\frac{G(x, y)}{g(x)g(y)}\asymp \frac{G(x_1, y)}{g(x_1)g(y)}\asymp \frac{\Psi(d(x,y))}{g(A)^2V(x, d(x,y))}.\end{equation} This establishes the $\mc{B}$-approximation property in case (f).\\

\noindent\textbf{(g)} For the final case, suppose that $2M \delta(y)<d(x,y)$. In particular, $r(x,y)=d(x,y)$. Take \begin{equation}r:=\frac{d(x,y)\wedge \epsilon}{2}, \textrm{   } x_1 \in \mc{A}_{r/M}(\xi), \textrm{    } y_1\in \mc{A}_{r/M}(\eta).\end{equation}

Then $x, x_1\in B(\xi, r/M)$ and $y, y_1\in B(\eta, r/M)$. Moreover, $d(y, y_1)\leq \delta(y)+r/M$ and $d(y, \xi)\geq d(x,y)-\delta(x)$ whence $$d(y_1, \xi)-r\geq d(y, \xi)-r-d(y, y_1)$$ $$\geq d(x,y)-r-\delta(x)-\delta(y)-\frac{r}{M}>\left(1-\frac{1}{2}-2\cdot \frac{1}{2M}-\frac{1}{2M}\right)d(x,y)\geq 0.$$ Similarly, \begin{equation}d(x_1, \eta)-r\geq d(x,\eta)-r-d(x,x_1)>0.\end{equation} By the Uniform Boundary Harnack Principle, \begin{equation}\frac{G(x,y)}{g(x)} \asymp \frac{G(x_1, y)}{g(x_1)}\end{equation} and \begin{equation}\frac{G(y, x_1)}{g(y)}\asymp \frac{G(y_1, x_1)}{g(y_1)}\end{equation} and therefore, by quasi-symmetry, \begin{equation}\frac{G(x,y)}{g(x)g(y)}\asymp \frac{G(x_1, y_1)}{g(x_1)g(y_1)}.\end{equation}

Moreover, \begin{equation}\abs{d(x_1, y_1)-d(x,y)}\leq d(x_1, x)+d(y_1,y)<\delta(x)+\delta(y)+\frac{2r}{M}<\frac{2}{M}d(x,y)\end{equation} whence $\frac{1}{3}d(x,y)<d(x_1, y_1)<2d(x,y)$. Obviously, $\epsilon d(x,y)\leq 2r$, since $d\leq 1$ and $\epsilon\leq 1$. Therefore $d(x_1, y_1)\leq 4r/\epsilon\leq (4M^2/\epsilon)\delta(x_1)\wedge \delta(y_1)$ and, by Lemma \ref{hansenl2.2} as well as (VD) and ($\Psi$D), \begin{equation}G(x_1, y_1)\asymp \frac{\Psi(d(x_1, y_1))}{V(x_1, d(x_1, y_1))}\asymp \frac{\Psi(d(x,y))}{V(x, d(x,y))}.\end{equation}

If $r=\epsilon/2$, i.e. if $d(x,y)\geq \epsilon$, then $g(x_1)\asymp g(y_1)\asymp g(A)\asymp 1$ by Lemma \ref{amazinglemma}. If $r<\epsilon/2$, then $d(x,y)<\epsilon$ and therefore $g(x_1)\asymp g(y_1)\asymp g(A)$ by the Harnack inequality \ref{bigassumption}(iv). Thus \begin{equation}\frac{G(x_1, y_1)}{g(x_1)g(y_1)}\asymp \frac{\Psi(d(x,y))}{V(x,d(x,y))g(A)^2}\end{equation} and the proof of Theorem \ref{mainpart2} is finished.\qed\\

Now it is time to prove Theorem \ref{mainpart3}. We begin with a lemma that is a generalized version of Carleson's estimate.\\

\begin{lemma}\label{hansenl2.6} Suppose that $D$ satisfies the Uniform Boundary Harnack Principle with respect to $G$. Let $g=1\wedge G(\cdot, o)$ be the benchmark Green's function. Then there exists a constant $C>0$ such that:

\noindent(i) For all $\xi\in \partial D$, every $0<r<\epsilon$, and every $A\in \mc{A}_{r/M}(\xi)$, \begin{equation}g\leq C g(A) \textrm{ on } D\cap B(\xi, r/M).\end{equation}

\noindent(ii) If $\xi \in \partial D$, $0<s\leq r<\epsilon$ and $A\in \mc{A}_r(\xi)$, then \begin{equation}g\leq C g(A) \textrm{ on } D\cap B(\xi, Ms)\cap \set{\delta>s/M}.\end{equation}

\noindent(iii) If $x,y,z\in D$ satisfy $d(z,x)\leq d(z,y)$, then \begin{equation}g(A)\leq Cg(B)\end{equation} for all $A\in \mc{B}(x,y)$ and $B\in \mc{B}(z,y)$.\end{lemma} 

\noindent\textbf{Proof:} (i) Fix $\xi \in \partial D$, $0<r<\epsilon$, $A\in \mc{A}_{r/M}(\xi)$, and choose $A'\in \mc{A}_{2Mr}(\xi)$. Then $r<d(A,A')<3Mr<3M^2\delta(A)\wedge \delta(A')$ and therefore, by Lemma \ref{hansenl2.2}, \begin{equation}G(A,A')\asymp \frac{\Psi(d(A,A'))}{V(A, d(A,A'))}.\end{equation}

For every $x\in D\cap B(\xi, r/M)$, $d(x, A')>r$ whence $G(x, A')\leq C_0 \frac{\Psi(d(x,A'))}{V(x,d(x,A'))} \asymp G(A,A')$, the latter because $d(x,A)\asymp r\asymp d(A,A')$ and the above application of Lemma \ref{hansenl2.2}.\\

The uniform boundary Harnack principle implies that \begin{equation}\frac{g(x)}{g(A)}\asymp\frac{G(x, o)}{G(A, o)}\leq C\frac{G(x, A')}{G(A,A')}\preceq C.\end{equation}%analyze this and fill in details (Laurent)

\noindent(ii) Fix $\xi \in \partial D$, $0<s\leq r<\epsilon$ and $A\in \mc{A}_r(\xi)$. Let $x\in D\cap B(\xi, Ms)$ such that $\delta(x)>s/M$. Choose $\widetilde{x}\in \mc{A}_{s/M}(\xi)$ and $\widetilde{A}\in \mc{A}_{r/M}(\xi)$. Then $g(\widetilde{x})\leq Cg(\widetilde{A})$ by (1). Moreover, $\delta(\widetilde{x})>s/M^2$ whence \begin{equation}d(x,\widetilde{x})\leq d(x,\xi)+d(\widetilde{x}, \xi)\leq (M+1)s\leq 2M^3\delta(x)\wedge \delta(\widetilde{x}).\end{equation} Thus $g(x)\asymp g(\widetilde{x})$ by Lemma \ref{hansenl2.2}. In particular, taking $r$ instead of $s$, $g(A)\asymp g(\widetilde{A})$.\\

\noindent(iii) Fix $x,y,z\in D$ such that $d(z,x)\leq d(z,y)$. Then $d(x,y)\leq d(x,z)+d(z,y)\leq 2d(z,y)$ and $\delta(x)\leq \delta(z)+d(z,x)\leq \delta(z)+d(z,y)$ and therefore \begin{equation}r(x,y)=\delta(x)\vee \delta(y)\vee d(x,y)\leq 2(\delta(z)\vee \delta(y)\vee d(z,y))=2r(z,y).\end{equation}\\

If $r(z,y)\geq \epsilon/2$, then $g\asymp 1$ on $\mc{B}(z,y)$ by Lemma \ref{amazinglemma}, and the conclusion holds since $g(A)\leq c_0$.\\

So suppose that $r(z,y)<\epsilon/2$ and therefore $r(x,y)<\epsilon$. Choose $\eta\in \partial D$ with $d(y,\eta)=\delta(y)$ and take $A\in \mc{A}_{r(x,y)}(\eta)$ and $B\in \mc{A}_{r(z,y)}(\eta)$. Using Lemma \ref{cones} we have \begin{equation}\mc{A}_{r(x,y)}(\xi)\cup\mc{A}_{r(x,y)}(\eta)\subseteq \mc{B}(x,y)\end{equation} and analogously for $(z,y)$. Therefore we have $A\in \mc{B}(x,y)$ and $B\in \mc{B}(z,y)$. Moreover, taking $s=r(x,y)/2$ we have $s\leq r(x,z)$, $\delta(A)>r(z,y)/M>s/M$ and $d(A, \eta)<r(x,y)\leq Ms$. Thus $g(A)\leq C g(B)$ by (2). As discussed above, the values of $g$ are comparable on any set of the form $\mc{B}(x,y)$. The conclusion follows.\qed\\

\noindent\textbf{Proof of Main Theorem \ref{maithm}, $\Leftarrow$ direction:} Assume that $D$ satisfies the Strong Generalized $3G$ Principle with respect to $G$. We wish to prove that $D$ satisfies the Uniform Boundary Harnack Principle with respect to $G$. Recall the alternate formulation of the Generalized $3G$ Principle given in Proposition \ref{alt3g}. Namely, that if $\widetilde{G}(x,y)=G(x,y)/(g(x)g(y))$, then we have that \begin{equation}\frac{\widetilde{G}(x,z)\widetilde{G}(y,z)}{\widetilde{G}(x,y)}\leq C\left(\widetilde{G}(x,z)+\widetilde{G}(z,y)\right)\end{equation} for all $x,y,z\in D$.\\

Fix $\xi\in \partial D$, $0<r\leq r_0$, and $x_1, x_2, y_1, y_2\in D$ with $Md(x_j, \xi)<r$ and $d(y_j, \xi)\geq r$, as in the setup for the Uniform Boundary Harnack Principle. Since $M\geq 3$, we have \begin{equation}d(x_1, x_2)<\frac{2}{3}r\leq d(x_1, y_1)\wedge d(x_2, y_2).\end{equation}

So the Strong Generalized $3G$ Principle implies that \begin{equation}\widetilde{G}(x_1, y_1)\leq C\widetilde{G}(x_2, y_1)\end{equation} and \begin{equation}\widetilde{G}(x_2, y_2)\leq C\widetilde{G}(x_1, y_2).\end{equation} Therefore \begin{equation}\frac{G(x_1, y_1)}{G(x_2, y_1)}\colon \frac{G(x_1, y_2)}{G(x_2, y_2)} = \frac{\widetilde{G}(x_1, y_1)\widetilde{G}(x_2, y_2)}{\widetilde{G}(x_2, y_1)\widetilde{G}(x_1, y_2)} \leq C^2.\end{equation} This finishes the $\Leftarrow$ direction of the Main Theorem \ref{maithm}.\qed\\

\textbf{Proof of Main Theorem \ref{maithm}, $\Rightarrow$ direction:} Assume that $D$ satisfies the Uniform Boundary Harnack Principle with respect to $G$. We wish to prove that $D$ satisfies the Strong Generalized $3G$ Principle with respect to $G$.\\

Fix $x,y,z\in D$ such that $d(z,x)\leq d(z,y)$ and therefore $d(x,y)\leq 2d(z,y)$. By Lemma \ref{hansenl2.6}, there exist $A\in \mc{B}(x,y)$ and $B\in \mc{B}(z,y)$ such that \begin{equation}g(A)\leq C g(B).\end{equation} By Theorem \ref{mainpart2}, $\mc{B}$-approximation holds. Thus \begin{equation}\frac{\widetilde{G}(z,y)}{\widetilde{G}(x,y)}\preceq  \frac{\Psi(d(z,y))}{g(B)^2 V(z, d(z,y))}\cdot \frac{g(A)^2V(x,d(x,y))}{\Psi(d(x,y))} \preceq 1.\end{equation} The last inequality holds due to Assumption \ref{bigassumption} (v). Notably this is our only use of (v) in the proof of the Main Theorem. \qed\\

To finish this section, we relate the Generalized $3G$ Principle to the quantity $\frac{\Psi(d(x,y))}{V(x,d(x,y))}$, which in future sections will be comparable to a global Green's function.\\

\begin{lemma}\label{global3g}
Assume that $D$ satisfies the Strong Generalized $3G$ Principle with respect to $G$. Then we have that for all $x,y,z\in D$, $$\frac{G(x,z)G(y,z)}{G(x,y)}\preceq \frac{\Psi(d(x,z))}{V(x,d(x,z))}+\frac{\Psi(d(y,z))}{V(y,d(y,z))}.$$
\end{lemma}

\noindent\textbf{Proof:} We begin by stating the weaker version, the Generalized $3G$ Principle: there is a constant $C>0$ such that for all $x,y,z\in D$, \begin{equation}
\frac{G(x,z)G(y,z)}{G(x,y)}\leq C\left(\frac{g(z)}{g(x)}G(x,z)+\frac{g(z)}{g(y)}G(z,y)\right).
\end{equation}

We will use this, as well as the fact that the Strong Generalized $3G$ Principle implies the $\mc{B}$-approximation Property and Carleson's Lemma. In light of the above, it suffices to prove that \begin{equation}
\frac{g(z)}{g(x)}G(x,z)\preceq \frac{\Psi(d(x,z))}{V(x,d(x,z))}.
\end{equation}

If $r(x,z)\geq \epsilon$, then by \ref{amazinglemma}, $g(z)\asymp g(x)\asymp 1$. The result then follows from Assumption \ref{bigassumption} (iii).\\

Therefore assume that $r(x,z)<\epsilon$. By the $\mc{B}$-approximation property, for any $A\in \mc{B}(x,y)$ we have $$
\frac{G(x,z)}{g(x)g(z)}\preceq \frac{\Psi(d(x,z))}{V(x,d(x,z))g(A)^2}$$

or $$
\frac{g(z)}{g(x)}G(x,z)\preceq \frac{\Psi(d(x,z))g(z)^2}{V(x,d(x,z))g(A)^2}.$$ To finish we must show that $g(z)\preceq g(A)$. Let $\xi\in \partial D$ be such that $d(z, \xi)=\delta(z)$. Note $\delta(z)\leq r(x,z)<\epsilon$. Set $r=\delta(z)$ and find $0<s\leq r$ such that $Ms>r$. By the Interior Corkscrew Property \ref{bigassumption}(i), there exists $A\in \mc{A}_r(\xi)$, which notably also satisfies $A\in \mc{B}(x,z)$ by Lemma \ref{cones}. Since $\delta(z)>r/M>s/M$, we can apply Lemma \ref{hansenl2.6}(ii) to conclude that $g(z)\preceq g(A)$ as desired.\qed\\

\section{Inner Uniform Domains and the Interior Corkscrew Property}
A central class of examples for the main theorem will occur on inner uniform domains. To define them, let us first start with length spaces.

\begin{definition}
Let $(X, d)$ be a metric space. The length $L(\gamma)\in[0,\infty]$ of a continuous curve $\fct{\gamma}{[0,1]}{X}$ is given by \begin{equation}
L(\gamma)=\sup\sum_i d(\gamma(t_{i-1}), \gamma(t_i)),
\end{equation} where the supremum is taken over all partitions $0=t_0<t_1<\cdots<t_k=1$ of $[0,1]$. $(X,d)$ is a \textbf{length space} if $d(x,y)$ is equal to the infimum of the lengths of continuous curves joining $x$ to $y$.
\end{definition}

\begin{definition}
A metric space $(X, d)$ is a \textbf{geodesic space} if for any $x,y\in X$ there is a continuous curve $\fct{\gamma}{[0,1]}{X}$ such that $\gamma(0)=x$, $\gamma(1)=y$, and $d(\gamma(s), \gamma(t))=\abs{s-t}d(x,y)$ for all $s, t\in [0,1]$.
\end{definition}

Next we introduce the intrinsic distance on a domain $D\subseteq X$.

\begin{definition}
Let $D\subseteq X$ be a domain. We define the \textbf{intrinsic distance} $d_i$ by \begin{equation} d_i(x,y)=\inf\set{L(\gamma)\colon \fct{\gamma}{[0,1]}{D} \textrm{ continuous, } \gamma(0)=x, \gamma(1)=y}.\end{equation}
\end{definition}

It is well known that $(D, d_i)$ is a length space. We let $(\widetilde{D}, d_i)$ denote its completion. For $x\in \widetilde{D}$ we let $B_{\widetilde{D}}(x,r)=\set{y\in \widetilde{D}\colon d_i(x,y)<r}$ and $B_D(x,r)=D\cap B_{\widetilde{D}}(x,r)$.\\

We let $\partial_{\widetilde{D}}D=\widetilde{D}\setminus D$, with accompanying distance to the boundary $\delta_{\widetilde{D}}(x)=d_i(x, \partial_{\widetilde{D}}D)$. If $V\subseteq D$ is open, let $\overline{V}^{d_i}$ denote the completion of $V$ with respect to $d_i$. We then write $$\partial_{\widetilde{D}} V=\overline{V}^{d_i}\setminus V.$$

\begin{definition}
Let $D$ be a domain in a complete, separable, locally compact length space. Let $\fct{\gamma}{[0,1]}{D}$ be a rectifiable, continuous curve in $D$. Let $c_D, C_D>0$. We say that $\gamma$ is a \textbf{$(c_D, C_D)$-inner uniform curve} if \begin{equation} L(\gamma)\leq C_Dd_i(\gamma(0), \gamma(1)),\end{equation} and \begin{equation}
\delta_D(\gamma(t))\geq c_D d_i(\gamma(0), \gamma(t))\wedge d_i(\gamma(t), \gamma(1))
\end{equation} for all $0\leq t\leq 1$. The domain $D$ is called an \textbf{inner uniform domain} if any two points in $D$ can be joined by a $(c_D, C_D)$-inner uniform curve.
\end{definition}

The following lemma yields the existence of inner uniform curves between any two points in $\widetilde{D}$.

\begin{lemma}
Let $(X, d)$ be a complete, locally compact, separable length space. Let $D$ be a $(c_D, C_D)$-inner uniform domain. Then for any $x,y\in \widetilde{D}$ there exists a $(c_D, C_D)$-inner uniform curve joining $x$ to $y$ in the $d_i$ metric.
\end{lemma}

\noindent\textbf{Proof:} See \cite{barlowmurugan1}, Lemma 2.6.\qed

\begin{example}
Let $D\subseteq \Real^n$ be bounded and Lipschitz. Then $D$ is inner uniform.
\end{example}

\begin{example}\label{slitdisk}
Let $D=\set{(x,y)\in \Real^2\colon x^2+y^2<1}\setminus\set{(x,0)\colon 0\leq x<1}$. Then $D$ is inner uniform.
\end{example}
 The above example \ref{slitdisk} is an important example of an inner uniform domain because it is not uniform, i.e. the inner metric cannot be replaced by the Euclidean metric in the definition of inner uniform curves.

\begin{proposition}\label{cork}
Let $D$ be an inner uniform domain in a complete, separable, locally compact length space. Then $D$ satisfies the Interior Corkscrew Condition as a subset of $\widetilde{D}$: namely, there exists $M\geq 3$ and $r_0>0$ such that for every $\xi\in \partial_{\widetilde{D}}D$ and all $0<r\leq r_0$, there exists $x\in D$ such that $d_i(x, \xi)<r$ and $\delta_D(x)>r/M$.
\end{proposition}

\noindent\textbf{Proof:} Fix $o\in D$ and set $r_0=\delta_D(o)/2$. Let $\xi\in \partial_{\widetilde{D}}D$. Let $\gamma$ be a $(c_D, C_D)$-inner uniform curve connecting $\xi$ to $o$. For $0<r\leq r_0$, by the intermediate value theorem there exists $t$ with $x=\gamma(t)$ such that $d_i(x,\xi)=r$. By inner uniformity of $\gamma$, we have $\delta_D(x)\geq c_D d_i(\xi, x)\wedge d_i(x, o)$. By the triangle inequality, $d_i(x,o)>r$. Thus $d_i(\xi, x)\wedge d_i(x,o)=r$. Setting $M=3\vee (c_D)^{-1}$ finishes the proof.\qed\\

Thus Assumption \ref{bigassumption}(i) is established for inner uniform domains. To discuss the classical Harnack principle, Assumption \ref{bigassumption}(iv), seems to by definition require a Green's function. We can circumvent this by phrasing the Harnack principle in terms of chains of balls; when the Green's function finally does come into play, the classical elliptic Harnack inequality in a ball will do the rest.

\section{Dirichlet Spaces and Applications}

In this section, we introduce the concept of Dirichlet spaces. See \cite{fukushima1} for a more complete discussion of Dirichlet spaces.\\% (Use Mathav's new paper \cite{murugan1})\\

\begin{definition}
A \textbf{metric measure Dirichlet (MMD) space} is a tuple $(X, d, \mu, \mc{E}, \mc{F})$, where $(X, d)$ is a locally compact, separable metric space, $\mu$ is a Radon measure on $X$ with full support, and $\mc{E}$ is a regular strictly local Dirichlet form on $L^2(X, d\mu)$ with domain $\mc{F}$. If $(\mc{E}, \mc{F})$ is local or strictly local we refer to the Dirichlet space as the same.\\

\begin{example}\label{classicalDirichlet}
The classical Dirichlet space is $\Real^n$ with the Euclidean metric and Lebesgue measure $dx$ with form \begin{equation}
\mc{E}(f,f)=\int_{\Real^n} \abs{\grad f}^2 dx
\end{equation}
with domain $\mc{F}=W^{1,2}(\Real^n, dx)$. This is a metric measure Dirichlet space.
\end{example}
\end{definition}

%Let $(X, d, \mu, \mc{E}, \mc{F})$ be a Dirichlet space. Let $\mc{L}$ be the generator of $(\mc{E}, \mc{F})$ in $L^2(X, d\mu)$ i.e. $\mc{L}$ is a self-adjoint, non-positive-definite operator with domain $\mc{D}(\mc{L})$ dense in $\mc{F}$ with \begin{equation}
%\mc{E}(f,g)=-\ip{\mc{L}f}{g}_{L^2}
%\end{equation}
%for all $f\in \mc{D}(\mc{L}))$ and all $g \in\mathcal{F}$. In Example \ref{classicalDirichlet}, $\mc{L}=\Delta$. The operator $\mc{L}$ has associated heat semigroup \begin{equation}
%P_t=e^{t\mc{L}}, \hspace{5mm} t\geq 0.
%\end{equation}

If the heat semigroup $(P_t)_{t>0}$ of $(\mc{E}, \mc{F})$ has an integral kernel $p(t,x,y)$ with respect to $\mu$ we call this the \textbf{heat kernel}. By definition this satisfies \begin{equation} P_tf(x)=\int_X f(y)p(t,x,y)d\mu(y)\end{equation} for all $f\in L^2(X, d\mu)$. In the classical case \ref{classicalDirichlet} we have \begin{equation}
p(t,x,y)=\frac{1}{(4\pi t)^{n/2}} \exp\left(-\frac{\norm{x-y}^2}{4t}\right).
\end{equation}

We mention the heat kernel because it is closely related to the Green's function for a Dirichlet space \begin{equation}
G(x,y):=\int_0^\infty p(t,x,y)dt
\end{equation} provided this integral is finite. Estimates for the heat kernel also yield estimates on the Green's function which will be essential in applying our Main Theorem \ref{maithm}.\\ 

To get the necessary bounds from Assumption \ref{bigassumption}(iii), we will impose bounds on the heat kernel of $X$. We will need to study scale functions and volume growth once again. Indeed, let $\fct{\Psi}{[0,\infty)}{[0, \infty)}$ be a scale function and define for $s>0$ \begin{equation}
\Phi(s)=\sup_{r>0} \left(\frac{s}{r}-\frac{1}{\Psi(r)}\right).
\end{equation}
 For example, if $\Psi(r)=r^\beta$, then $\Phi(s)=s^{\frac{\beta}{\beta-1}}\left(\frac{1}{\beta^{\frac{1}{\beta-1}}}-\frac{1}{\beta^{\frac{\beta}{\beta-1}}}\right).$

\begin{definition}
We say that an MMD space $(X,d,\mu, \mc{E}, \mc{F})$ satisfies (full) \textbf{sub-Gaussian heat kernel estimates (HKE($\Psi$))} with respect to $\Psi$ if $\Psi$ is a scale function and there are constants $c_i>0$ such that for all $t>0$ and all $x,y\in X$ the heat kernel $p(t,x,y)$ of $X$ satisfies $$
\frac{c_1}{V(x, \Psi^{-1}(t))}\exp\left(-c_2t\Phi\left(c_3\frac{d(x,y)}{t}\right)\right)\leq p(t,x,y),$$ $$p(t,x,y)\leq \frac{c_4}{V(x, \Psi^{-1}(t))}\exp\left(-c_5t\Phi\left(c_6\frac{d(x,y)}{t}\right)\right).
$$
\end{definition}

The heat kernel estimates HKE($\Psi$) have profound implications. We collect the ones that will be useful to us here.

\begin{definition}
$(X, d, \mu, \mc{E}, \mc{F})$ satisfies the \textbf{(scale-invariant) elliptic Harnack inequality (EHI)} if there exist $C>0$ and $0<\theta<1$ such that for any open ball $B$ in $X$ and any $u\in \mc{F}$ such that $u\geq 0$ and $u$ is harmonic in $B$ we have $$\esssup_{\theta B} u\leq \essinf_{\theta B} u.$$ 
\end{definition}

Note that our harmonic functions will have continuous versions, so in practice we refer to maxima and minima instead of essential suprema and infima.

\begin{definition}\label{pipsi}
$(X, d, \mu, \mc{E}, \mc{F})$ satisfies the \textbf{Poincar\'e inequality PI($\Psi$)} on $X$ if $\Psi$ is a scale function and there exist constants $C_P>0$ and $a_P\geq 1$ such that for any ball $B$ in $X$ and any $f\in \mc{F}$, $$\int_B\abs{f-f_B}^2d\mu\leq C_P\Psi(R)\int_{a_P B} d\Gamma(f, f).$$
Here $f_B=\mu(B)^{-1}\int_B fd\mu$ and $\Gamma(f,f)$ is the energy measure of $f$.
\end{definition}

\begin{definition}\label{respsi}
$(X, d, \mu, \mc{E}, \mc{F})$ satisfies the \textbf{resistance condition RES($\Psi$)} if there exists a constant $C_{res}>1$ such that for any ball $B$ with radius $R>0$ such that $2B\neq X$, $$C_{res}^{-1}\frac{\mu(B)}{\Psi(R)}\leq \textrm{Cap}_{2B}(B)\leq C_{res}\frac{\mu(B)}{\Psi(R)}.$$ Here $\textrm{Cap}$ is capacity.
\end{definition}

\begin{proposition}
Let $(X, d, \mu, \mc{E}, \mc{F})$ be an MMD space that satisfies HKE($\Psi$). Then $(X, d, \mu, \mc{E}, \mc{F})$ satisfies volume doubling, EHI, PI($\Psi$), and RES($\Psi$).
\end{proposition}

\noindent\textbf{Proof:} See Lierl \cite{lierl1}, Theorem 2.12.\qed\\

For our purposes, HKE($\Psi$) alone is not enough to capture the desired behavior of Green's functions, and so we must impose further requirements on $\Psi$ and $V$.\\

\begin{assumption}\label{bigassumption2}
Assume that there exists $\gamma>0$, $C>0$, and $A>0$ such that for our scale function $\Psi$ and volume profile $V$ and all $0<r\leq R\leq \textrm{diam}(X)/A$ we have \begin{equation}
C^{-1}\left(\frac{R}{r}\right)^2\leq \frac{\Psi(R)}{\Psi(r)}\leq C\left(\frac{r}{R}\right)^{\gamma} \frac{V(R)}{V(r)}.
\end{equation}
\end{assumption}

In $\Real^n$ for $n\geq 3$, the scale function $\Psi(r)=r^2$ and volume profile $V(r)=r^n$ satisfies Assumption \ref{bigassumption2} with $\gamma=n-2$. In $\Real$ and $\Real^2$, \ref{bigassumption2} fails. This assumption may therefore be viewed as a form of transience or non-parabolicity.\\

The significance of this assumption for us is twofold. First, in a Dirichlet space satisfying HKE($\Psi$), we can integrate the heat kernel to obtain an estimate on the Green's function of the space, if that integral is finite. Second, there exist examples of spaces satisfying HKE($\Psi$) for a wide range of scale functions and volume profiles, thus giving new examples of our Main Theorem.\\

\subsection{Proof that Assumption \ref{bigassumption} holds in inner uniform domains under HKE($\Psi$) and Assumption \ref{bigassumption2}}

In this section we demonstrate that our main assumption \ref{bigassumption} holds for a large array of inner uniform domains in MMD spaces. We begin with a connectness result essential to proving the classical Harnack inequality \ref{bigassumption}(iv). This is similar to Proposition 2.8 in \cite{minerbe} which concerns manifolds; we have adapted this result to a Dirichlet space setting.

\begin{lemma}\label{rca}
Let $(X, d, \mu, \mc{E}, \mc{F})$ be an MMD space satisfying HKE($\Psi$) and Assumption \ref{bigassumption2}. Assume also that $(X, d)$ is a length space. Then there exists $\kappa\geq 2$ and $A>0$ such that for any $o\in X$ and $0<R<\textrm{diam}(X)/A$, the annulus $B(o, R)\setminus \overline{B}(o, 2^{-1}R)$ is contained in a connected component of $B(o, R)\setminus \overline{B}(o, \kappa^{-1}R)$.
\end{lemma}

\noindent\textbf{Proof:} We will use that HKE($\Psi$) implies volume doubling as well as PI($\Psi$) (\ref{pipsi}): there exist $C_P>0$ and $a_P\geq 1$ such that for all $x\in X$, $r>0$, and all $f\in \mc{F}$, 

\begin{equation}
\int_{B(x,r)} \abs{f-f_{B(x,r)}}^2d\mu\leq C_P\Psi(r)\int_{B(x, a_Pr)} d\Gamma(f,f).
\end{equation}
Here $f_B=\mu(B)^{-1}\int_B fd\mu$ and $\Gamma(f,f)$ is the energy measure of $f$.\\

Let $A$ be such that for any $0<R<\textrm{diam}(X)/A$ and any $o\in X$, the annulus $B(o, R)\setminus \overline{B}(o, R/2)$ is nonempty. Let $0<R<\textrm{diam}(X)/A$ be fixed. For $i\geq 0$ set $A_i=B(o, 2^{-i}R)\setminus \overline{B}(o, 2^{-(i+1)}R)$. Set $B_i=B(o, 2^{-i}R)$. If $A_0$ is connected, we are done with $\kappa=2$. Otherwise let $X$ and $Y$ be connected components of $A_0$. Let $I$ be the least integer, if it exists, such that $X$ and $Y$ are contained in the same connected component of $B(o, R)\setminus \overline{B}(o, 2^{-I-1}R)$. If $2^{-I+3}\geq a_P^{-1}$ we are once again done with $\kappa=2^{I+2}$, so henceforth we can assume that $2^{-I+3}<a_P^{-1}$. If no such $I$ exists, we let $I$ be arbitrary.\\

By definition, $X$ and $Y$ are in different connected components of $B(o, R)\setminus \overline{B}(o, 2^{-I}R)$. let $X'$ and $Y'$ denote these components and let $Z=X'\cap B(o, 2^{-I+1}R)$.\\

By RES($\Psi$) (\ref{respsi}), which is a consequence of HKE($\Psi$), there exists $\phi \in \mc{F}$ such that $\phi\in \mc{F}$, $\phi=1$ on $B(o, 2^{-I+1}R)$, $\phi=0$ outside of $B(o, 2^{-I+2}R)$, and \begin{equation}\int_{B(o, 2^{-I+2}R)\setminus B(o, 2^{-I+1}R)}d\Gamma(\phi, \phi)=\mc{E}(\phi,\phi)\leq C_{res}\frac{V(2^{I}R)}{\Psi(2^{I}R)}.\end{equation} The constant $C_{res}$ is independent of $o$, $I$ and $R$.\\

We now define $f$ as follows: $$
f=\begin{cases} 
          1 & \textrm{ on }X'\setminus Z \\
          1-\phi & \textrm{ on }Z \\
	0 & \textrm{ else}\\
       \end{cases}.$$

Multiplying $f$ by a cutoff function that is $1$ in $B(o, 2^{-1}R)$ and $0$ outside $B(o, R)$ if necessary, we may assume that $f\in \mc{F}$. Applying the Poincar\'e inequality to $f$, we see that \begin{equation}
\int_{B(o, 2^{-1}a_P^{-1} R)}\abs{f-f_{B(o, 2^{-1}a_P^{-1}R)}}^2d\mu\leq C_P\Psi(2^{-1}a_P^{-1}R)\int_{B(o, 2^{-1}R)}d\Gamma(f,f).
\end{equation}

Note that $2^{-I+4}< 2^{-1}a_P^{-1}$. Let $J\in \N$ be least such that $2^{-J}<2^{-1}a_P^{-1}$. Let $X_J=X'\cap A_J$ and $Y_J=Y'\cap A_J$. By connectedness, there exist points $x_J\in X_J$ and $y_J\in Y_J$ such that $d(o, x_J)=d(o, y_J)=3\cdot 2^{-J-2}R$. By considering balls $B(x_J, 2^{-J-3}R)$ and $B(y_J, 2^{-J-3}R)$ and applying volume doubling, we find that there is a constant $C_{VD}>0$ such that $V(o, 2^{-1}a_P^{-1}R)\leq C_{VD}(\mu(X_J)\wedge \mu(Y_J))$.\\

Furthermore, since $J\leq I-4$, we have that $f=1$ on $X_J$ and $f=0$ on $Y_J$. Therefore setting $B=B(o, 2^{-1}a_P^{-1} R)$ we have $$\int_{B}\abs{f-f_{B}}^2d\mu\geq \frac{\int_{B}\int_{B}\abs{f(x)-f(y)}^2d\mu(x)d\mu(y)}{2\mu(B)}\geq \frac{\mu(X_J)\mu(Y_J)}{2\mu(B)}\geq \left(\frac{1}{2C_{VD}^2}\right)\mu(B).$$

Note that $f$ is locally constant on $B(o, 2^{-1}R)\setminus Z$ so using the strict locality of the Dirichlet form $\mc{E}$, $$\int_{B(o, 2^{-1}R)}d\Gamma(f, f)=\int_{Z}d\Gamma(f, f)\leq\mc{E}(\phi, \phi)\leq C_{\res}\frac{V(2^{-I}R)}{\Psi(2^{-I}R)}.$$

The Poincar\'e inequality gives us that $$\left(\frac{1}{2C_{VD}^2}\right)\mu(B)\leq C_P\Psi(2^{-1}a_P^{-1}R)C_{\res}\frac{V(2^{-I}R)}{\Psi(2^{-I}R)}.$$ Moving it all to one side gives $$1\leq 2C_{VD}^2C_PC_{res}\left(\frac{\Psi(2^{-1}a_P^{-1}R)}{\Psi(2^{-I}R)}\right)\left(\frac{V(2^{-I}R)}{V(2^{-1}a_P^{-1}R)}\right).$$ Applying \ref{bigassumption2} gives $$1\leq 2C_{VD}^2C_PC_{res}\left(\frac{2^{-I}}{2^{-1}a_P^{-1}}\right)^{\gamma}.$$ Since $\gamma>0$ and all constants involved are independent of $o$, $R$, and $I$, we must have that $I$ is bounded above. This also proves that such an $I$ must exist, as if it does not, we can use this argument for arbitrarily large $I$ to get a contradiction.\qed\\

The conclusion of the above lemma is sometimes known as the relatively connected annuli (RCA) condition. We now build the necessary concepts to apply the elliptic Harnack inequality.

\begin{definition}\label{hchain} $(X, d)$ be a metric space and let $U\subseteq X$ be a connected open set and $M\geq 1$. For $x,y\in U$, an \textbf{$M$-Harnack chain} from $x$ to $y$ in $U$ is a sequence of balls $B_1, \dots, B_n$ each contained in $U$ such that $x\in M^{-1}B_1, y\in M^{-1}B_n$, and $M^{-1}B_j\cap M^{-1}B_{j+1}\neq \emptyset$ for each $j=1, \dots, n-1$. The number $n$ of balls in a Harnack chain is called its length. Write $N_U(x,y;M)$ for the length of the shortest $M$-Harnack chain in $U$ from $x$ to $y$.
\end{definition}

\begin{lemma}\label{chainlemma}
Let $(X, d, \mu, \mc{E}, \mc{F})$ be an MMD space satisfying HKE($\Psi$). Let $D\subseteq X$ be a $(c_D, C_D)$-inner uniform domain such that $\partial D\neq \emptyset$. Then for each $M>1$ there exists $C>0$ depending only on $c_D, C_D$ and $M$ such that for all $x,y\in U$ \begin{equation}
C^{-1}\log\left(\frac{d_i(x,y)}{\delta_D(x)\wedge \delta_D(y)}+1\right)\leq N_D(x,y;M)\leq C\log\left(\frac{d_i(x,y)}{\delta_D(x)\wedge \delta_D(y)}+1\right) +C.
\end{equation}
\end{lemma}

\noindent\textbf{Proof:} See Barlow and Murugan \cite{barlowmurugan1}, Lemma 4.7.\qed\\

\begin{lemma}\label{annh}
Let $(X, d, \mu, \mc{E}, \mc{F})$ be an MMD space satisfying HKE($\Psi$) and Assumption \ref{bigassumption2} and such that $(X, d)$ is a length space. There exists a constant $C>0$ such that if $o\in X$ and $R>0$ and $u>0$ is a positive harmonic function in $X\setminus \set{o}$, then $u(x)\leq Cu(y)$ for any $x,y\in B(o, R)\setminus \overline{B}(o, R/2)$.
\end{lemma}

\noindent\textbf{Proof:}  Let $\kappa>0$ be as in Lemma \ref{rca}. We now demonstrate the existence of $N\in \N$ such that for any $o\in X$ and $R>0$, the annulus $A=B(o,  R)\setminus \overline{B}(o, \kappa^{-1}R)$ can be covered by $N$ balls of radius at most $4^{-1}\kappa^{-1}R$.\\

Indeed, let $B_1, \dots, B_M$ be a maximal pairwise disjoint set of balls of radius $16^{-1}\kappa^{-1}\theta R$ that meet $A$. Here $\theta$ is as in the statement of the elliptic Harnack inequality. By maximality, the balls $5B_1, \dots, 5B_M$ cover $A$. Therefore $\mu(A)\leq \sum_{j=1}^M \mu(5B_j)\leq C\sum_{j=1}^M \mu(B_j)$. However, by the triangle inequality, the balls $B_j$ are all contained in the ball $B(o, (1+8^{-1}\kappa^{-1})R)$ which has volume $\asymp B(o,  R)$ by volume doubling. Thus we have that $M\mu(B(o, 16^{-1}\kappa^{-1} R)\preceq\sum_{j=1}^M \mu(B_j)\leq \mu(A)\preceq B(o,  R)$. We then conclude by volume doubling that $M\preceq 1$ with the constant involved depending only on that given from volume doubling. This proves the existence of such an $N$ claimed above.\\

So let $B_1, \dots, B_N$ be balls of radius $16^{-1}\kappa^{-1}\theta R$ as in the claim. Let $x,y\in B(o, r)\setminus \overline{B}(o, r/2)$. Using Lemma \ref{rca}, we can build a $\theta$-Harnack chain from $x$ to $y$ in $X\setminus \set{o}$ of balls from $\set{B_1, \dots, B_N}$ that naturally has length $\leq N$. Note that $\theta^{-1}B_j\subseteq X\setminus \set{o}$ for each $j$, so the elliptic Harnack inequality (EHI) applies to each ball $B_j$ with a uniform constant. The conclusion of the lemma follows.\qed\\

We now state the result that allows us to apply the main theorem to inner uniform domains.

\begin{theorem}\label{appthm}
Let $(X, d, \mu, \mc{E}, \mc{F})$ be an MMD space such that $(X, d)$ is a geodesic space and let $\Psi$, $V$ be a scale function and volume profile for $X$ satisfying Assumption \ref{bigassumption2}. Assume also that HKE($\Psi$) holds. Let $D$ be a bounded inner uniform domain in $X$ such that $\textrm{diam}(D)< \textrm{diam}(X)/4$ and let $G$ be its Green's function. Then Assumption \ref{bigassumption} holds with respect to $D$ and $G$.
\end{theorem}

\noindent\textbf{Remark:} The condition that $\textrm{diam}(D)< \textrm{diam}(X)/4$ can be weakened considerably; we choose it here for simplicity of proof and of use. For example, one could replace this with the condition that the first Dirichlet eigenvalue of $D$ within $X$ is strictly positive and use results from Lierl \cite{lierl1}.\\

\noindent\textbf{Proof:} (i) The Interior Corkscrew Property holds by Proposition \ref{cork}.\\

(ii) Quasi-symmetry is trivial: the Green's function is symmetric.\\

(iii)  We start with the lower bound. First let $x, y\in D$ be distinct such that $d_i(x,y)=d(x,y)\leq \delta_d(y)/2$. Assume first that $d(x,y)>(3/7)\delta_D(y)$. Let $\Gamma$ be a path from $x$ to $y$ of length $d(x,y)$ and let $z$ lie on the image of $\Gamma$ with $d(z, y)=\delta_D(y)/7$. Let $R=(3/7)\delta_D(y)$ and note that $B(y, 2R)\neq X$. By Lemma 3.13 in Lierl \cite{lierl1}, the Green's function $G_{B(y, R)}$ of $B(y, R)$ satisfies $G_{B(y, R)}(y,z)\geq C\frac{\Psi(d(y,z))}{V(z, d(y,z))}$, where the constant $C$ is independent of $z$ and $y$. If $d(x,y)<R$ then by the same as well as domain monotonicity of Green's functions we have $G(x,y)\geq G_{B(y, R)}(x,y)\geq C\frac{\Psi(d(x,y))}{V(x, d(x,y))}$ as desired. Using a Harnack chain of balls centered on $\Gamma$, we can then use EHI to conclude that $G(x,y)\preceq G_D(z,y) \preceq C\frac{\Psi(d(x,y))}{V(x, d(x,y))}$ as desired. Finally if $d(x,y)\leq (3/7)\delta_D(y)$ then no such chain is necessary and we repeat the argument with $z=x$.\\

Now for the upper bound. Let $x,y\in D$ be distinct. Since $\textrm{diam}(D)<\textrm{diam}(X)/4$, if $R=\textrm{diam}(D)$ then $B(x, 2R)\neq X$ and $y\in B(x,R)$. The upper bound in Lemma 3.13 in Lierl \cite{lierl1} then yields exactly the result we desire.\\

(iv) Let $y\in D$, let $k\in \N$ and let $x_1, x_2\in D(y)$ such that $d(x_1, x_2)\leq k \delta(x_1)\wedge \delta(x_2)$. By Lemma \ref{chainlemma}, there is a $\delta^{-1}$-Harnack chain $B_1, \dots, B_N$ joining $x_1$ to $x_2$ of length at most $C\log(k+1)+C$.\\

Case 1: $y\notin B_j$ for all $j$.\\

Here we can apply the scale-invariant elliptic Harnack inequality EHI in each ball $\delta^{-1}B_j$ to conclude that $G(x_1, y)\leq C^{C\log(k+1)+C} G(x_2, y)$ as desired.\\

Case 2: $y\in B_j$ for some $j$.\\

Let $R=\delta_D(y)/2\kappa$, $\kappa>0$ from Lemma \ref{rca}, and let $A=B(y, R)\setminus B(y, \kappa^{-1}R)$. Let $S=\set{x\colon d(x,y)=R}$. $j_1=\min\set{j\colon S\cap B_j\neq \emptyset}$ and $j_2=\max\set{j\colon S\cap B_j\neq \emptyset}$. Note that $j_1$ and $j_2$ exist by connectness. Write $B_{j_1}=B(x_{j_1}, r_{j_1})$ and $B_{j_2}=B(x_{j_2}, r_{j_2})$.\\

If $y\notin B_{j_1}\cup B_{j_2}$ we can apply Lemma \ref{annh} to get that $\sup_{B_{j_1}}G(\cdot, y)\leq \inf_{B_{j_2}} G(\cdot, y)$, and finish the proof with applying the elliptic Harnack inequality across the Harnack chains from $B_1$ to $B_{j_1}$ and from $B_{j_2}$ to $B_N$.\\

If $y\in B_{j_1}$, we note by the triangle inequality that $r_1\leq 2\delta_D(y)$. Furthermore, since $y\in B_{j_1}$ and $S\cap B_{j_1}\neq \emptyset$, we conclude that $r_1\geq \frac{1}{4}\delta_D(y)$. By a covering argument similar to the proof of Lemma \ref{annh}, there exists a universal constant $Q$ such that we can cover $\delta^{-1}B_{j_1}$ in $Q$ balls $B$ of radius $\frac{1}{5\delta}\delta_D(y)$. Note that for each such ball $B$, we have $\delta B\subseteq D$ and $y\notin \delta B$.\\

We can therefore , via a simple recursion argument, find a $\delta^{-1}$-Harnack chain of at most $Q$ balls connecting $B_{j_1-1}$ to $S$ while avoiding $y$. Therefore we have an $\delta^{-1}$-Harnack chain of length $j_1-1+Q$ from $x_1$ to $S$ which avoids $y$. A similar argument yields an $\delta^{-1}$-Harnack chain of length $N-j_2+Q$ from $S$ to $x_2$ which avoids $y$. Applying the elliptic Harnack inequality to these chains as well as Lemma \ref{annh} to $S$ yields the Harnack inequality $G(x_1, y)\leq C(k)G(x_2, y)$ as desired.\qed\\

(v) This follows from Assumption \ref{bigassumption2}.\qed\\

Thus, if $D$ is a domain in $(X, d, \mu, \mc{E}, \mc{F})$ under the hypotheses of Theorem \ref{appthm}, by our Main Theorem \ref{maithm}, the Uniform Boundary Harnack Principle is equivalent to the Strong Generalized $3G$ Principle on $D$. We can further state that both of these are simply true, due to the following theorem of Barlow and Murugan \cite{barlowmurugan1}.

\begin{theorem}\label{bestboundaryharnack} (Barlow, Murugan 2019) Let $D$ be an inner uniform domain in $(X, d, \mu, \mc{E}, \mc{F})$ under all the same hypotheses on $D$ and $(X, d, \mu, \mc{E}, \mc{F})$ as Theorem \ref{appthm}. Let $G(x,y)$ be the Green's function of $D$. Then $D$ satisfies the Uniform Boundary Harnack Principle with respect to $G$.
\end{theorem}

We therefore have the following immediate corollary.
\begin{corollary}
Let $D$ be an inner uniform domain in $(X, d, \mu, \mc{E}, \mc{F})$ under all the same hypotheses on $D$ and $(X, d, \mu, \mc{E}, \mc{F})$ as Theorem \ref{appthm}. Let $G(x,y)$ be the Green's function of $D$. Then $D$ satisfies the Strong Generalized $3G$ Principle with respect to $G$.

\end{corollary}
\noindent\textbf{Proof:} This follows from combining our Main Theorem \ref{maithm} with Theorem \ref{appthm} as well as the above Theorem \ref{bestboundaryharnack}.\qed

\subsection{Examples: Fractals and Generalized Fractals}

\begin{figure}
\centering
\includegraphics[scale=0.1]{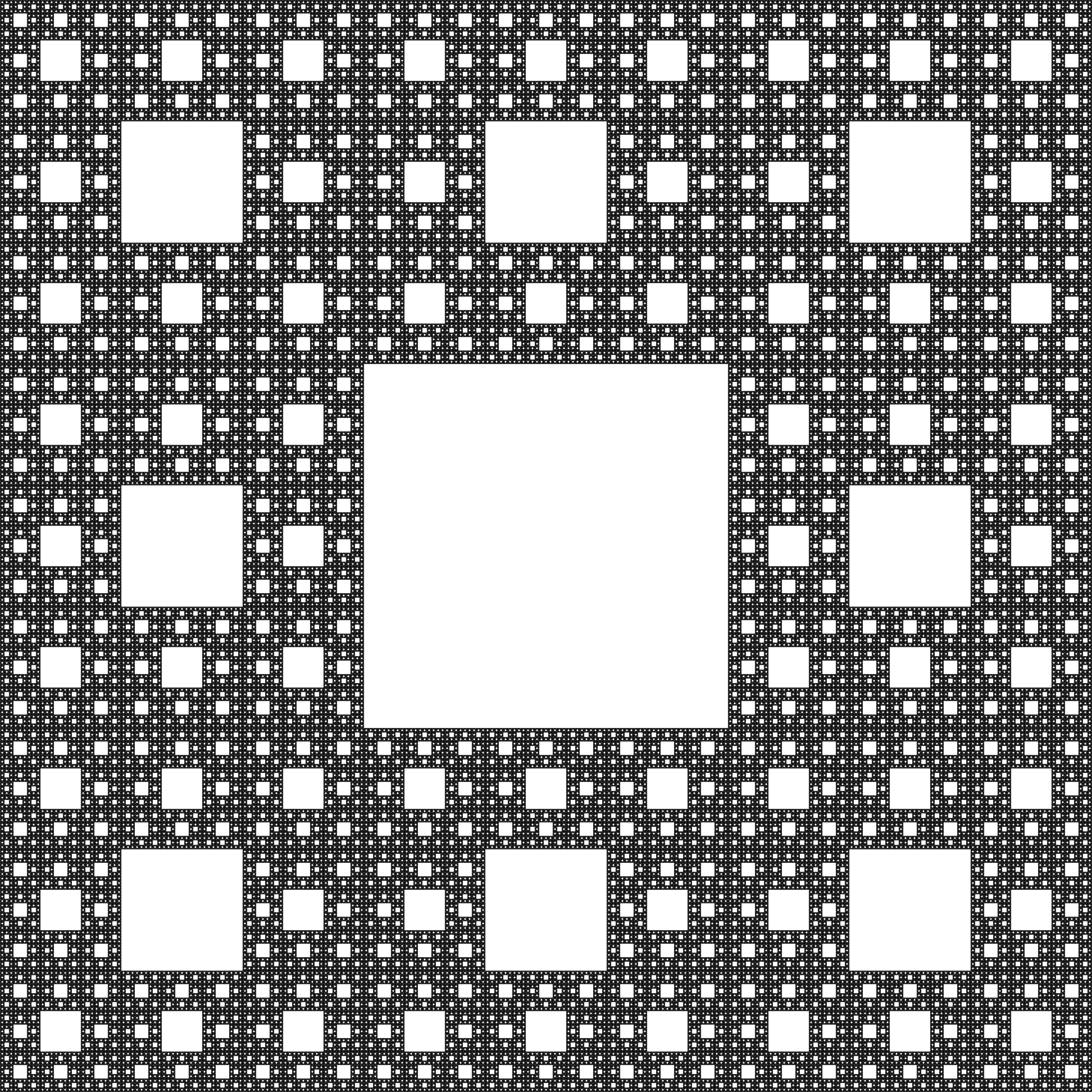}% no need to specify the file extension
\caption{Compact Sierpinski carpet in $\Real^2$. Public domain image from Wikimedia Commons.}
\label{fig:2dcarpet}
\end{figure}

Let us consider the Sierpinski carpet. Our presentation is taken from \cite{barlow1}. Let $L=3$, $I=3^n-1$ and consider $[0,1]^n$ as a union of the $3^n$ closed cubes of the form $\prod_{j=1}^n[a_j/3, (a_j+1)/3]$ for $a_j=0,1,2$. Enumerate these cubes $C_0, C_1, \dots, C_{3^{n}-1}$ with $C_0=[1/3, 2/3]^n$. Let $\psi_i,$ $i=1, \dots, I$ be the orientation-preserving affine maps that take $F_0=[0, 1]^n$ onto $C_i$. Set $\psi_1(x)=L^{-1}x$. Define $\mc{P}(A)=\bigcup_{i=1}^I \psi_i(A)$ for $A\subseteq \Real^n$. Taking an $m$-fold composition we define $F_m=\mc{P}^{(m)}(F_0)$ and then $F=\bigcap_{m=0}^\infty F_m$. This yields the compact Sierpinski carpet, as pictured in Figure \ref{fig:2dcarpet} for dimension $n=2$. Notably we will focus on higher dimensions as $n=2$ fails Assumption \ref{bigassumption2}. $F$ has Hausdorff dimension $d_f=\log(3^n-1)/\log 3$; let $\mu_F$ denote $d_f$-dimensional Hausdorff measure on $F$. $F$ is a length space with length metric $d$ comparable to the Euclidean metric. $(F, d, \mu_F)$ is a metric measure space.\\

For a more comprehensive study of Dirichlet forms on the Sierpinski carpet, see \cite{barlowbass1}, \cite{barlowbass2}, \cite{barlowbasskumagaiteplyaev1}, and \cite{barlow1}. There are still many open questions, but notably by \cite{barlowbasskumagaiteplyaev1}, there exists, up to a constant, a unique regular strictly local Dirichlet form on $F$ that is invariant under the symmetries of $F$. Let $\mc{E}$ denote this Dirichlet form and let $\mc{D}$ denote its domain. Notably, there exists $\beta>2$ such that $(F, d, \mu_F, \mc{E}, \mc{D})$ satisfies HKE($\Psi$) with $\Psi(r)=r^\beta$. This quantity $\beta$ is known as the \textbf{walk dimension} of $F$, due to it being almost surely the Hausdorff dimension of the trace of an interval of Brownian motion in $F$ provided $\beta\leq d_f$. Furthermore, if $n\geq 3$, then indeed $\beta<d_f$ and therefore Assumption \ref{bigassumption2} holds.\\

There are many examples of inner uniform domains for which our results apply. For one of particular note, let $E=\set{x=(x_1, \dots, x_n)\in [0, 1]^n\colon x_j=1\textrm{ for some }j=1, \dots, n}$. Then $E\subseteq \partial [0,1]$ and $E\subseteq F$. Consider the domain $D=F\setminus E$ in $F$ with boundary $E$. This notably fails our condition on diameter, since $D$ has the same diameter as $F$. However, $(1/9)D=\set{(1/9)x\colon x\in D}$ is a domain in $F$ with sufficiently small diameter, and hence by our results, the Strong Generalized $3G$ Theorem holds. But by the fractal symmetry of $F$, the Dirichlet form on $(1/9)D$ is a scalar multiple to that on $D$ itself, so $D$ satisfies the Strong Generalized $3G$ Theorem.\\

It is notable to compare our result with the paper of Hansen \cite{hansenmain}. The main result of that paper can be used in place of our main Theorem; however, our Theorem \ref{appthm} is still required to show that the assumptions of Hansen's result hold on the Sierpinski carpet.\\

For our next example, we turn to the paper of Murugan \cite{murugan1}, which constructs MMD spaces with prescribed heat kernel estimates.

\begin{theorem}\label{prescribed} (Murugan, 2024) Let $\fct{V, \Psi}{[0,\infty)}{[0, \infty)}$ be doubling functions and let $C>1$ be such that for all $0<r\leq R$,\begin{equation}\label{scalerel}
C^{-1}\frac{R^2}{r^2}\leq \frac{\Psi(R)}{\Psi(r)}\leq C\frac{RV(R)}{rV(r)}.
\end{equation} Then there exists an unbounded MMD space $(X, d, \mu, \mc{E}, \mc{F})$ that satisfies HKE($\Psi$) and such that $(X, d)$ is a geodesic space and $V$ is a volume profile for $(X, d, \mu, \mc{E}, \mc{F})$.

\end{theorem}

In particular, if instead $\frac{\Psi(R)}{\Psi(r)}\leq C\left(\frac{r}{R}\right)^\gamma \frac{V(R)}{V(r)}$ for some $\gamma>0$ and all $0<r\leq R$, then the resulting unbounded MMD space will satisfy Assumption \ref{bigassumption2}, and thus any bounded inner uniform domain in that space will be subject to the Main Theorem \ref{maithm}. We refer to Murugan's paper for the details of the construction, which is somewhat complicated. The spaces themselves are quotients of trees; trees by themselves have too little volume to satisfy Assumption \ref{bigassumption2}, but by gluing them together carefully one can control the volume profile in an adequate manner.\\

\begin{proposition}
Let $(X, d, \mu, \mc{E}, \mc{F})$ be an MMD space constructed via the proof of Theorem \ref{prescribed} in \cite{murugan1}. Then there exist points $x\in X$ and $c, C>0$ such that for all $k\in \Z$, the ball $B(x, 2^k)$ is $(c, C)$-inner uniform.
\end{proposition}

\noindent\textbf{Proof:} See \cite{murugan1}, Proposition 3.26.\qed\\

What is noteworthy about the spaces in Theorem \ref{prescribed} is that neither $\Psi(r)$ nor $V(r)$ need be strictly single powers of $r$. Rather they can be irregular, growing at different rates depending on the scale of $r$. Thus these spaces lack the symmetry of spaces such as the Sierpinski carpet, which has a well-defined Hausdorff dimension and well-defined walk dimension. See \cite{hansenmain} for a result similar to our Main Theorem, but which does not cover such irregular spaces as these.

\subsection{Application: Conditional Gaugability of Schr\"odinger Operators}

Let $D$ be a bounded inner uniform domain in an unbounded geodesic MMD space $(X, d, \mu, \mc{E}, \mc{F})$ satisfying HKE($\Psi$) and Assumption \ref{bigassumption2}. Let $G(x,y)$ be the Green's function of $X$ and let $G_D(x,y)$ be the Green's function of $D$. Since $X$ is unbounded, $G$ necessarily exists by Assumption \ref{bigassumption2}. In fact, by similar estimates to those appearing in the proof of Lemma 3.13 in \cite{lierl2}, we have that $G(x,y)\asymp \frac{\Psi(d(x,y))}{V(d(x,y))}$ globally. Recall Lemma \ref{global3g}, a form of the $3G$ Principle relating $G_D$ to $G$.\\

 Let $(B_t)_{t>0}$ be the Markov process on $X$ associated with $(\mc{E}, \mc{F})$. In the classical setting, this is Brownian motion, so we can think of $(B_t)_{t>0}$ as generalized Brownian motion. This process has transition kernel $p(t,x,y)$, where $p$ is the heat kernel. In $D$, we can consider the same process with killing at the boundary of $D$. This process has transition kernel $p_{Dir}(t,x,y)$, where $p_{Dir}$ is the Dirichlet heat kernel in $D$. We can further consider the conditional process in $D$ with transition kernel $G_D(x,y)^{-1}p_{Dir}(t,x,z)G_D(z,y)$. The Markov process associated with this transition kernel is $(B_t)_{t>0}$ in $D$ conditioned to end at $y$. We denote expectation with respect to this process started at $x$ and ending at $y$ by $E^x_y$. Let $\tau_D$ denote the lifetime of the process.\\

Let $\fct{W}{D}{\Real}$ be such that $\sup_{x\in X} \int_D G(x,y)\abs{W(y)}d\mu(y)<+\infty$. Set $\norm{W}=\sup_{x\in X} \int_D G(x,y)\abs{W(y)}d\mu(y)$. We can think of $W$ as the potential for a Schr\"odinger operator $-\mc{L}+W$, where $-\mc{L}$ is the generator of the heat semigroup on $X$.\\

Utilizing Lemma \ref{global3g}, we have that $$E^x_y\left[\int_0^{\tau_D} \abs{W(B_s)}ds\right]=\int_D \frac{G_D(x,z)G_D(z,y)}{G_D(x,y)}\abs{W(z)}d\mu(z)$$ $$\preceq \int_D(G(x,z)+G(z,y))\abs{W(z)}d\mu(z)\preceq \norm{W}.$$ Thus we have: if $\norm{W}$ is sufficiently small, then $$\sup_{x,y\in D}E^x_y\left[\int_0^{\tau_D} \abs{W(B_s)}ds\right]<1.$$ From here we can apply Khasminskii's Lemma to conclude that $W$ is conditionally gaugeable in $D$, i.e. $$\sup_{x,y\in D}E^x_y\left[\exp\left(\int_0^{\tau_D} \abs{W(B_s)}ds\right)\right]<+\infty.$$ This implies, for example, the existence of a Green's function for the Dirichlet space $(X, d, \mu, \mc{E}+W\ip{\cdot}{\cdot}, \mc{F})$. We refer the reader to Chen \cite{zqchen1} for more on the gaugeability of Schr\"odinger operators, as well as Hansen \cite{hansen2} for extensive applications of the $3G$ Principle to Schr\"odinger operators.

\subsection{Acknowledgements}

The first author was supported in part by the National Science Foundation Graduate Research Fellowship grants number DGE - 2139899 and DGE - 1650441. The second author was supported in part by National Science Foundation grants DMS-2054593 and DMS-2343868.

\nocite{aikawa3}

\bibliographystyle{plain}

\bibliography{3GPaperFinalVersionBibliography.bib}

\end{document}